\RequirePackage{ifpdf}
\ifpdf % We are running pdfTeX in pdf mode
\documentclass[pdftex]{sigma}
\else
\documentclass{sigma}
\fi

\newtheorem{mytheorem}{Theorem}[section]
\newtheorem{myprop}[mytheorem]{Proposition}
\newtheorem{mylemma}[mytheorem]{Lemma}
\newtheorem{mycor}[mytheorem]{Corollary}

\newcommand{\bg}{{\mathfrak g }}
 \newcommand{\ba}{{\mathfrak a }}

 \newcommand{\bl}{{\mathfrak l }}
 \newcommand{\bk}{{\mathfrak k }}
 \newcommand{\bq}{{\mathfrak q }}
 \newcommand{\bp}{{\mathfrak p}}
 \newcommand{\bu}{{\mathfrak u }}
 \newcommand{\bh}{{\mathfrak h }}
 \newcommand{\bt}{{\mathfrak t }}
 \newcommand{\bs}{{\mathfrak s}}

 \newcommand{\bR}{{\mathbb R}}
 \newcommand{\bN}{{\mathbb N}}
 \newcommand{\bC}{{\mathbb C}}
 \newcommand{\bZ}{{\mathbb Z}}

\newcommand{\bgo}{{\mathfrak g}_o}

\newcommand{\blo}{{\mathfrak l}_o}
\newcommand{\bko}{{\mathfrak k}_o}

\newcommand{\bho}{{\mathfrak h}_o}

\newcommand{\bso}{{\mathfrak s}_o}

\numberwithin{equation}{section}

\begin{document}

%\allowdisplaybreaks

\renewcommand{\PaperNumber}{017}

\FirstPageHeading

\renewcommand{\thefootnote}{$\star$}

 \ShortArticleName{Branching Laws for Some Unitary Representations of
 $SL(4,{\mathbb R})$}

 \ArticleName{Branching Laws for Some Unitary Representations\\ of
 $\boldsymbol{SL(4,{\mathbb R})}$\footnote{This paper is a
contribution to the Proceedings of the 2007 Midwest
Geometry Conference in honor of Thomas~P.\ Branson. The full collection is available at
\href{http://www.emis.de/journals/SIGMA/MGC2007.html}{http://www.emis.de/journals/SIGMA/MGC2007.html}}}

 \Author{Bent \O{}RSTED~$^\dag$ and Birgit SPEH~$^\ddag$}

 \AuthorNameForHeading{B.~\O{}rsted and B.~Speh}

 \Address{$^\dag$~Department of Mathematics, University of Aarhus, Aarhus,
Denmark}
 \EmailD{\href{mailto:orsted@imf.au.dk}{orsted@imf.au.dk}}
 \URLaddressD{\url{http://person.au.dk/en/orsted@imf}}

 \Address{$^\ddag$~Department of Mathematics,  310 Malott Hall,
 Cornell University,\\
 $\phantom{^\ddag}$~Ithaca, NY 14853-4201,  USA}
 \EmailD{\href{mailto:speh@math.cornell.edu}{speh@math.cornell.edu}}

\ArticleDates{Received September 10, 2007, in f\/inal form January
27, 2008; Published online February 07, 2008}

 \Abstract{In this paper  we consider the restriction of a unitary
 irreducible representation of type $A_{\bq}(\lambda)$ of
 $GL(4,\bR)$ to reductive subgroups $H$ which are the f\/ixpoint sets
 of an involution.  We obtain a formula for the restriction to the
 symplectic group and to $GL(2,\bC)$, and as an application we
 construct in the last section some representations in the cuspidal
 spectrum of the symplectic and the complex general linear group.
 In addition to working directly with the cohmologically induced
 module to obtain the branching law, we also introduce the
 useful concept of pseudo dual pairs of subgroups in a reductive
 Lie group.}

 \Keywords{semisimple Lie groups; unitary representation; branching
 laws}

 \Classification{22E47; 11F70}

\section{Introduction}

Understanding a unitary representation $\pi$ of a Lie groups $G$
often involves understanding its restriction to suitable subgroups
$H$. This is in physics referred to as breaking the symmetry, and
often means exhibiting a nice basis of the representation space of
$\pi$. Similarly, decomposing a tensor product of two
representations of $G$ is also an important branching problem,
namely the restriction to the diagonal in $G \times G$. Generally
speaking, the more branching laws we know for a given
representation, the more we know the structure of this
representation. For example, when $G$ is semisimple and $K$ a
maximal compact subgroup, knowing the $K$-spectrum, i.e.\ the
collection of $K$-types and their multiplicities, of $\pi$ is an
important invariant which serves to describe a good deal of its
structure. It is also important to give good models of both $\pi$
and its explicit $K$-types. There has been much progress in recent
years (and of course a large number of more classical works, see
for example \cite{W,hs,howe,JV}), both for
abstract theory as in \cite{K, ko1,ko3,ko2}, and concrete examples of branching laws in \cite{var,varhol,loke,varfa,zhang}.

In this paper, we shall study in a special case a generalization
of the method applied in~\cite{martens} and again in~\cite{JV};
this is a method of Taylor expansion of sections of a vector
bundle along directions normal to a submanifold. This works nicely
when the original representation is a~holomorphic discrete series
for $G$, and the subgroup $H$ also admits holomorphic discrete
series and is embedded in a suitable way in $G$. The branching law
is a discrete sum decomposition, even with f\/inite multiplicities,
so-called admissibility of the restriction to $H$; and the
summands are themselves holomorphic discrete series
representations for $H$. Since holomorphic discrete series
representations are cohomologically induced representations in
degree zero, it is natural to attempt a generalization to other
unitary representations of similar type, namely cohomologically
induced representations in higher degree. We shall focus on the
line bundle case, i.e.\ the $A_{\bq}(\lambda)$ representations. In
this case T.~Kobayashi~\cite{ko2} obtained necessary and
suf\/f\/icient conditions  that the restriction is discrete and that
each representation appears with f\/inite multiplicity,
socalled {\it admissibility} of the representation
relatively to the subgroup. Using
explicit resolutions and f\/iltrations associated with the imbedding
of $H$ in $G$, we analyze the derived functor modules and obtain
an explicit decomposition into irreducible representations. It is
perhaps not surprising, that with the appropriate conditions on
the imbedding of the subgroup, the class of (in our case derived
functor) modules is preserved in the restriction from $H$ to $G$.

While the algebraic methods of derived functor modules, in particular
the cohomologically induced representations, provide a very strong
tool for the theoretical investigations of unitary representations
of reductive Lie groups, it has been dif\/f\/icult to work with
concrete models of these modules. It is however exactly these
models that we use in this paper, as outlined above; the fact
that one may consider these modules as
Taylor expansions
of appropriate dif\/ferential forms, indicates that it should be
natural to study these Taylor expansions along submanifolds~--
and when these submanifolds are natural for the subgroup for
which one wants to do the branching law, there arises a useful
link between the algebraic branching and the geometry of the
imbedding of the subgroup. It is our hope, that this idea (that
we carry out in some relatively small examples) will have a
broader use in deciding the possible candidates for representations
occuring in an admissible branching law.

Here is the general setting that we consider: Let $G$ be a
semisimple linear connected Lie group with maximal compact
subgroup $K$ and Cartan involution $\theta $. Suppose that
$\sigma$ is another involution so that $\sigma \cdot \theta =
\theta \cdot \sigma$ and let $H$ be the f\/ixpoint set  of $\sigma $
in $G$. Suppose that $L = L_x$  is  the centralizer of an elliptic
element $x \in i(\bg \cap \bh )$ and  let $\bq = \bl \oplus \bu$,
$\bq^H = \bq \cap \bh $ be the corresponding $\theta$-stable
parabolic subgroups. Here we use as usual gothic letters for
complex Lie algebras and subspaces thereof; a subscript will
denote the real form, e.g.~$\bgo$. We say that pairs of parabolic
subalgebras $\bq$, $\bq^H$  which are constructed this way are well
aligned. For a unitary character~$\lambda $ of~$L$ we def\/ine
following Vogan/Zuckerman the unitary representations
$A_{\bq}(\lambda)$.

In this paper we consider the example of the group $G=SL(4,\bR)$.
There are two $G$-conjugacy classes of skew symmetric matrices
with representants $Q_1=\left(
\begin{array}{cc} J&0\\ 0 & -J \end{array}\right)$    and
$Q_2=\left(
\begin{array}{cc} J & 0\\ 0 & J \end{array}\right)$ where $J=\left(
\begin{array}{cc} 0 & -1 \\ 1& 0 \end{array}\right)$. Let $H_1$
respectively $H_2$ be the  symplectic subgroups def\/ined by these
matrices and $H_1'$, $H_2'$ the centralizer of $Q_1$, respectively
$Q_2$. All these subgroups are f\/ix point sets of involutions
$\sigma_i$, $i = 1,2$ and $\sigma_i'$, $i = 1,2$ respectively.

The matrix $Q_2$ has f\/inite order, it is contained in all
subgroups $H_i$ and $i Q_2 \in \bg$ def\/ines a~$\theta $-stable
parabolic subalgebra $\bq$ of~$sl(4,\bC$) and also
$\theta$-stable parabolic subalgebras $\bq_{\bh_1 }=\bq \cap
\bh_1 $ of~$\bh_1$, respectively $\bq_{\bh_2}= \bq \cap \bh_2 $ of
$\bh_2 $. Its centralizer $L$ in $SL(4,\bR)$ is isomorphic to
$GL_1(2,\bC)= \{T\in GL(2,\bC)| \ |\det(T)| =1\}.$ The parabolic
subgroups $\bq$, $\bq_{\bh_1}$ as well as $\bq$, $\bq_{\bh_2}$ are
well aligned.

\looseness=1
We consider in this paper the unitary representation $A_{\bq}$ of
$G$ corresponding to trivial cha\-racter~$\lambda$. Its
inf\/initesimal character is the same as that of the trivial
representation. The representation~$A_\bq$ was studied from an
analytic point of view by S.~Sahi \cite{sahi}. Since the $A_\bq$
has nontrivial $(\bg,K)$-cohomology and is isomorphic to a
representation in the residual spectrum, this representation is
also interesting from the point of view of automorphic forms. See
for example~\cite{speh}. We show in this paper that the
restriction of $A_\bq$ to $H_1$ and $H_1'$ is a direct sum of
irreducible unitary representations, where as the restriction to
$H_2$ and $H_2'$ has continuous spectrum. We also determine
explicitly  the restriction of $A_\bq $ to the subgroups $H_1$ and
$ H_1'$ and conclude that for all unitary $(\bh_1, H_1\cap K)$-modules $V$ the dimension of Hom$_{(\bh_1,H_1\cap K)}( A_\bq,V)$
is at most 1.

If we interpret $SL(4,\bR)$ and $Sp(2,\bR)$ as Spin(3,3) and
Spin(2,3) then these branching laws can in some sense be
considered as supporting the conjectures by B. Gross and D. Prasad
\cite{gross-prasad} for the restriction of Vogan packets of
representations of $SO(n,n)$ to representations of $SO(n-1,n)$.

 The paper is organized as follows:
After introducing all the notation in Section~\ref{sec1}  we prove in
Section~\ref{sec2} using a result of T.~Kobayashi, that the restriction of
$A_{\bq}$ to $H_1$ and $H_1'$ is a direct sum of irreducible
unitary representations, whereas the restriction to $H_2$ and
$H_2'$ does have a~continuous spectrum. This discrete/continuous
alternative, see~\cite{ko1}, is one of the deep results that we
invoke for symmetric subgroups.

We do not attempt in this paper to say anything about the
continuous spectrum, and we mainly focus on the admissible
situation, so the alternative is really admissible/non-admissible.

In Sections~\ref{sec3} and~\ref{sec4} we determine the representations of $H_1$
respectively of $H_1'$ that appear in the restriction of $A_{\bq}$ to $H_1$ respectively $H_1'$ and show that it is a direct
sum of unitary representations of the form $A_{\bq \cap \bh _1
}(\mu )$ respectively $A_{\bq \cap \bh _1' }(\mu' )$, each
appearing with multiplicity one. The main point is here, that we
f\/ind a natural model in which to do the branching law, based on
the existence results of T.~Kobayashi; and also following
experience from some of his examples, where indeed derived functor
modules decompose as derived functor modules (for the smaller
group).

In Section~\ref{sec5} we introduce pseudo dual pairs. This allows us to
f\/ind another interpretation of the restrictions of $A_\bq$ to the
pseudo dual pair $H_1$, $H_2$.
This notion turns out to be extremely useful for analyzing the
spectrum in the admissible situation, and combined with our idea
of restricting a cohomologically induced module gives the
complete branching law. In Section~\ref{sec6} we recall some more examples of
                     branching laws.

In  Section~\ref{sec7} we formulate a conjecture about the multiplicity of
representations in the restriction of representations $A_{\bq }$
of semisimple Lie groups $G$ to subgroups $H$, which are
centralizers of involutions. If the restriction of $A_{\bq}$ to
$H$ is a direct sum of irreducible representation of $H$ we expect
that there is a $\theta$-stable parabolic subalgebra $\bq^H$ of
$H$ so that all representations which appear in the restriction
are of the form $A_{\bq^H}(\mu )$ and that a Blattner-type formula
holds. See the precise conjecture at the end of Section~\ref{sec7}, where
we introduce a natural generalization of previously known
Blattner-type formulas for the maximal compact subgroup.

In Section~\ref{sec8} we see how these results may be used to construct automorphic
representations of
 $Sp(2,\bR)$ and $GL(2,\bC)$ which are in the discrete spectrum
 for some congruence subgroup. For $Sp(2,\bR)$  these representations
 are in the residual spectrum, whereas for
  $GL(2,\bC)$ these
 representations are in the cuspidal spectrum.
We expect that our methods extend to other situations with similar
applications to automorphic representations; and we hope the point
of view introduced here also will help to understand in a more explicit
way the branching laws for semisimple Lie groups with respect to
reductive subgroups.

\section{Notation and generalities}\label{sec1}

\noindent {\bf 2.1.}  Let $G$ be a connected linear semisimple
Lie group. We f\/ix a maximal compact subgroup~$K$ and Cartan
involution~$\theta$. Let $H$ be a $\theta $-stable connected
semisimple subgroup with maximal compact subgroup $K^H=K \cap H$.
We pick a fundamental Cartan subgroup $C^H =T^H \cdot A^H$ of~$H.$
It is contained in a fundamental Cartan subgroup $C=T\cdot A$ of G
so that $T^H = T \cap H $ and $A^H = A \cap H.$ The complex Lie
algebra of a Lie group (as before) is denoted by small letters and
its real Lie algebra by a subscript $o.$ We denote the Cartan
decomposition by $\bg_o= \bk_o \oplus \bp$.

\medskip \noindent
{\bf Def\/inition.} Let $\bq$ and $\bq^H$ be $\theta$-stable
parabolic subalgebras of $\bg$, respectively $\bh $. We say that
they are {\it well aligned } if $\bq^H =\bq \cap \bh$.
\medskip

We f\/ix  $x_o$ in $\bt^H$. Then $i\ x_o$ def\/ines well aligned
$\theta $-stable parabolic subalgebras $\bq = \bl \oplus \bu$ and
$\bq ^H = \bl^H \oplus \bu^H =\bq \cap \bh$ of $\bg $ respectively
$\bh $; for details see page 274 in \cite{K-V}.

We write $L$ and $L^H$ for the centralizer of $x_0$ in $G$ and in
$H$ respectively. For a unitary character $\lambda$  of $L$ we
write $\lambda ^H$ for the restriction of $\lambda$ to ${L}^H$.

\medskip

\noindent {\bf 2.2.}  For later reference we recall the
construction of the representations $A_{\bq}(V)$, $V$ an
irreducible $(\bq, L\cap K) $ module. We follow conventions of
the book by Knapp and Vogan~\cite{K-V} (where much more detail on
these derived functor modules is to be found~-- this is our
standard reference) and will always consider representations of
$L$ and not of the metaplectic cover of $L$ as some other authors.
We consider $U(\bg) $ as right $U(\bq )$ module and write
$V^\sharp = V \otimes \wedge^{{\rm top}}\bu.$ Let $\bp_L$ be a~$L
\cap K$-invariant complement of $\bl \cap \bk$ in $\bl$. We write
$r_G= \bp_L \oplus \bu $. Now we introduce the derived functor
modules as on page~167 in~\cite{K-V},  recalling that this
formalizes Taylor expansions of certain dif\/ferential forms. Since
all the groups considered in the paper are connected we use the
original def\/inition of the Zuckerman functor~\cite{V1} and do not
use the the Hecke algebra $R(\bg,K)$ to def\/ine the representations
$A_{\bq}(V)$.
 Consider the complex
\begin{gather*}
0 \rightarrow \mbox{Hom}_{L \cap K
}(U(\bg),\mbox{Hom}(\wedge^0 r_G,V
^\sharp))_K \rightarrow \\
\phantom{0}{} \rightarrow
\mbox{Hom}_{L \cap K}(U(\bg),\mbox{Hom}(\wedge^1 r_G,V
^\sharp ))_K \rightarrow
% \\ \phantom{0}{} \rightarrow
 \mbox{Hom}_{L \cap K}(U(\bg), \mbox{Hom}(\wedge^2 r_G, V
^\sharp) )_K \rightarrow \cdots.
\end{gather*}
Here the subscript $K$ denotes the subspace of $K$-f\/inite vectors.
 We denote by $ T(x,U(  \cdot   ))$  an element in
 $\mbox{Hom}_{L\cap K}(U(\bg),\mbox{Hom}_{\bC}(\wedge^{n-1} r_G,V ^\sharp)
)_K$.
 The dif\/ferential $d$ is def\/ined by
 \begin{gather*}
 d\ T(x,U(X_1 \wedge X_2\wedge \cdots \wedge X_n))
 =\sum_{i=1}^{n} (-1)^i
 T(X_i x,U( X_1 \wedge X_2 \wedge \cdots \hat{X_i}\cdots \wedge X_n))\\
 \qquad {} + \sum_{i=1}^{n} (-1)^{i+1}
 T( x,X_iU( X_1 \wedge X_2 \wedge \cdots\hat{X_i}\cdots \wedge X_n))
\\
 \qquad {} + \sum_{i<j} (-1)^{i+j}
 T( x,U(P_{r_G}[X_i,X_j]\wedge X_1 \wedge X_2 \wedge \cdots \hat{X_i}\cdots
 \hat{X_j}
\cdots\wedge  X_n )),
\end{gather*}
 where $x\in U(\bg),\ X_j \in r_G$ and $P_{r_G}$ is the projection
 onto $r_G$ along $\bl\cap \bk$.
 Let $s= \dim (\bu \cap \bk)$ and let $\chi$
 be the inf\/initesimal character of $V$. If
 \[ \frac{2\langle \chi + \rho(\bu),\alpha\rangle}{|\alpha|^2} \not \in \{0,-1,-2,-3, \dots \}
 \qquad \mbox{for} \quad \alpha \in \Delta(\bu),\]
 where $\langle \ ,\ \rangle$ denotes the Killing form of $\bg$, then the cohomology is zero
 except in degree $s$ and if~$V$ is irreducible this def\/ines an
  irreducible $((U(\bg ),K)$-module
 $A_{\bq}(V)$ in degree~$s$ (8.28 in~\cite{K-V}).  By (5.24 in~\cite{K-V},
see also the remark/example on page~344) the inf\/initesimal
character of $A_{\bq}(V)$ is $\chi +\rho(\bu )$ (usual shift of
the half-sum of all positive roots in $\bu$).

If $V$ is trivial the inf\/initesimal character of $A_{\bq}(V)$ is
$\rho_G $ and we write simply
 $A_{\bq}$. Two representations $A_{\bq}$ and $A_{\bq'}$ are
equivalent if $\bq$ and $\bq'$ are conjugate under the compact
Weyl group $W_K$.

 For an irreducible
 f\/inite dimensional $(\bq ^H, L^H\cap K)$-module $V^{L ^H}$
 we  def\/ine similarly the  $(U({\bh} ),K^H)$-modules $A_{\bq ^H}(V^{L ^
H}) $.

\medskip

\noindent{\bf 2.3.} Let $H$ be the f\/ix point set  of an
involutive automorphism $\sigma$ of $G$ which commutes with the
Cartan involution $\theta$. We write $\bgo =\bho \oplus \bso $ for
the induced decomposition of the Lie algebra.
 T.~Kobayashi proved~\cite{ko1} that the
restriction of $A_{\bq}$ to $H$ decomposes as direct sum of
irreducible representations of $H$ if $A_\bq$ is $K^H$-admissible,
i.e.\ if every $K^H$-type has f\/inite multiplicity. If $A_{\bq}$ is
discretely decomposable as an $(\bho, K\cap H)$-module we call an
irreducible ($\bho, H\cap K$)-module~$\pi^H$ an {\it H-type} of
$A_{\bq}$ if {\samepage
\[\mbox{Hom}_{(\bho,K^H)} (\pi^H, A_{\bq}) \not = 0\]
  and the dimension of
 Hom$_{(\bho ,K^H)} (\pi^H, A_{\bq})$ its {\it multiplicity.} }

We have  $\bl =\bl^H \oplus \bl \cap \bs$. Put $\bu^H = \bu \cap
\bh$. The representation of $\bl^H$ on $\bu$ is  reducible and as
$\bl^H$-module $\bu= \bu^H \oplus (\bu \cap \bs)$. Let
$\overline{\bq }= \bl \oplus \overline{\bu} $ be the opposite
parabolic subgroup. Then $ \bh =\bk^H \oplus  \bu ^H\oplus
\overline{\bf u}^H$ and $ (\bu \cap \bs) \oplus \overline{\bf u}
\cap \bs$ is a $\bl^H$-module. As an $\bl^H$-module $\bg = \bh
\oplus (\bu \cap \bs) \oplus (\bl \cap \bs) \oplus
(\overline{\bu}\cap \bs)$.

\medskip

 \noindent{\bf 2.4.}  Now let $G= SL(4,\bR).$ The skew symmetric
matrices
\[Q_1=\left(
\begin{array}{cc} J & 0 \\ 0 & -J \end{array}\right)\qquad \mbox{and} \qquad Q_2=\left(
\begin{array}{cc} J & 0\\ 0 & J \end{array}\right)\] with $J=\left(
\begin{array}{cc} 0 & -1 \\ 1& 0 \end{array}\right)$ represent
the conjugacy classes of skew symmetric matrices
 under G.
 They def\/ine symplectic forms also denoted by
$Q_1$ and $Q_2$.

Let $H_1$, respectively $H_2$, be the $\theta $-stable symplectic
subgroups def\/ined by $Q_1$, respectively~$Q_2$. These subgroups
are f\/ix points of the involutions
\[\sigma_i(g) = Q_i\cdot (g^{-1})^{\rm tr} \cdot Q_i^{-1},\qquad i=1,2.
\]  Since $Q_1$ and $Q_2$ are conjugate in $GL(4,\bR )$, but not  in $SL(4,\bR)$,
 the symplectic groups $H_1$ and~$H_2$ are not conjugate in $SL(4,\bR)$.

 Let $H_1'$ and $H_2'$ be the f\/ix points of the involutions
\[\sigma_{i}'(g) = Q_i\cdot g\cdot  Q_i^{-1} .\]
Both groups $H_1'$ and $H_2'$ are  isomorphic to \[GL_1(2,\bC)=
\{T\in GL(2,\bC)|\  |\det(T)| =1 \}, \] but they are  not conjugate
in $SL(4,\bR)$.

\medskip

\noindent{\bf 2.5.}  We f\/ix $x_0=Q_2 $. It has f\/inite order and
is contained in  $\bigcap_{i=1}^{2}H_i$ and in
$\bigcap_{i=1}^{2}H_i'$.  Now $i x_0 \in i\bg$ def\/ines a $\theta $
stable parabolic subalgebra $\bq$ of $sl(4,\bC)$ and also
$\theta$-stable well aligned parabolic subalgebras $\bq^H$ of the
subalgebras $\bh$.
 Its centralizer $L= L_{x_0}$ in
$SL(4,\bR)$, the Levi subgroup, is isomorphic to $GL(2,\bC)
=H_2'$.  For a precise description of the parabolic see page~586
in~\cite{K-V}.

We have
\begin{gather*}
L=H_2',\qquad
K^{H_2} =  K\cap H_2 =K \cap L, \qquad
K^{H_2'}= K\cap H_2'  =K \cap L
\end{gather*}
and
\[ K^{H_1} = K^{H_1'} . \]

Let $ A_{\bq}$ be the representation holomorphically induced from
$\bq$ which has a  trivial inf\/ini\-te\-simal character. This
representations is a subrepresentation of a degenerate series
representation induced from a one dimensional representation of
the parabolic subgroup with Levi factor $S(GL(2,\bR)\times
GL(2,\bR))$ and thus all its $K$-types have multiplicity one. See
\cite{sahi} for details.

The next proposition demonstrates how dif\/ferent imbeddings of the
same subgroup (symplectic res. general linear complex) gives
radically dif\/ferent branching laws.

\begin{myprop}\label{propositionI.1}

\noindent
\begin{enumerate}\itemsep=0pt
 \item[\rm 1.] The restriction of $A_\bq $ to $H_1$ and to $H_1'$ is a
 direct sum of irreducible representation each
appearing with finite multiplicity.

\item[\rm 2.] The restriction of $A_\bq $ to  $H_2 $ and to $H_2'$
is not admissible and has
continuous spectrum.

\end{enumerate}
\end{myprop}

\begin{proof}
 Since $K^{H_1}=K^{H_1'}$ and $K^{H_2}=K^{H_2'}$
it suf\/f\/ices by T.~Kobayashi's Theorem~4.2 in~\cite{ko1} to  show that
$A_\bq$ is $K^{H_1}$ admissible. We will prove this in the next
section.

To prove (2) it suf\/f\/ices to show by  Theorem~4.2 in~\cite{ko1} that
$A_\bq $ is not $K^{H_2}$ admissible. This is proved also in the
next section.
\end{proof}

\noindent
{\bf Remark.}
We will prove later in the paper that the
representations of $H_1$ have at most  multiplicity one in the
restriction of $A_\bq $.

\section[The restriction of $ A_{\bq}$ to $K\cap H_i$, $i=1,2$]{The restriction of $\boldsymbol{A_{\bq}}$ to $\boldsymbol{K\cap H_i}$, $\boldsymbol{i=1,2}$}\label{sec2}

 \noindent We use in this section the notation
introduced on pages 586--588 in~\cite{K-V}.

\medskip

\noindent {\bf 3.1.}  The Cartan algebra $\bt_o$ of so($4,\bR)$
consists of 2 by 2 blocks $\left(
\begin{array}{cc} 0 & \theta_j\\ -\theta_j & 0
\end{array}\right)$ down the diagonal. We have a $\theta
$-stable Cartan subalgebra $ \bh_o=\bt_o \oplus \ba_o$ where
$\ba_o$ consists of the 2 by 2 blocks $\left(
\begin{array}{cc} x_j & 0\\ 0 & x_j \end{array}\right)$,
also down the diagonal.
We def\/ine $e_j \in \bh^*$ by
\[e_j\left(
\begin{array}{cc} x_j & -iy_j\\ iy_j & x_j \end{array}\right)=
y_j\] and $f_j \in \bh^*$ by
\[f_j\left(
\begin{array}{cc} x_j & -iy_j\\ iy_j & x_j
\end{array}\right)=x_j.\]
Then the roots $\Delta(\bu )$ of $(\bh,\bu)$ are
\[e_1 + e_2 +(f_1-f_2), \qquad e_1 + e_2 -(f_1-f_2), \qquad 2e_1, \qquad 2 e_2
\]
and a compatible set of positive roots $\Delta^+(\bl )$ of
$(\bh,\bl )$ are
\[e_1-e_2 +(f_1-f_2), \qquad e_1-e_2 -(f_1-f_2). \]
The roots $\alpha_1 = e_1 +e_2$, $\alpha_2 = e_1 -e_2$ are
compatible positive roots of the Lie algebra $\bk $ with respect
to $\bt$.

The highest weight of the minimal $K$-type of $A_\bq $ is
$\Lambda = 3( e_1 +e_2)$. See page 588 in~\cite{K-V}.  All other
$K$-types are of the form \[\Lambda +m_1(e_1 +e_2)+ 2m_2  e_1, \qquad  m_1, m_2 \in \bN.
\]

\begin{figure}[t]
\centering
\begin{minipage}[b]{75mm}
\centerline{\includegraphics{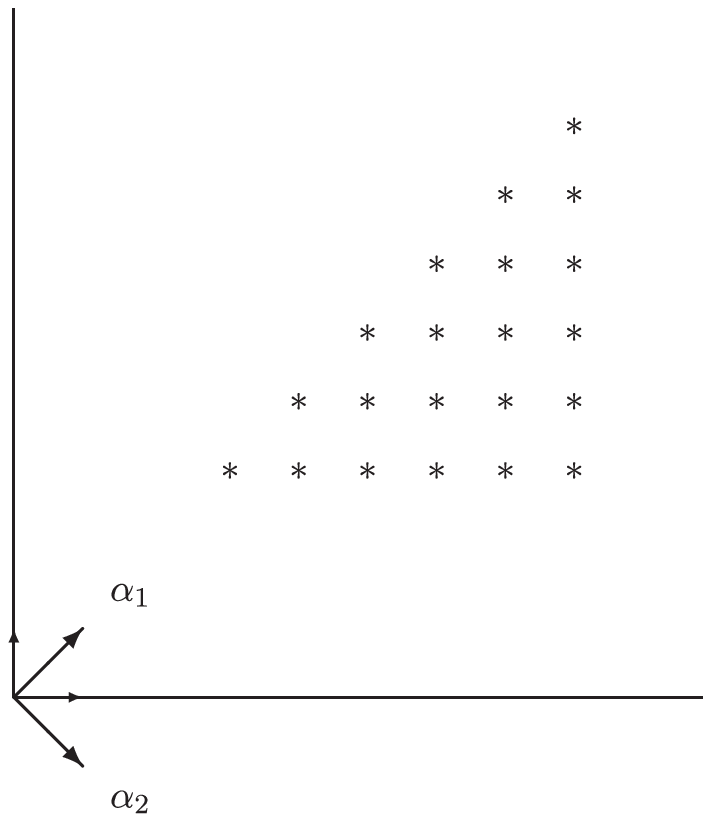}}
\caption{}\label{fig1}
\end{minipage}
\qquad
\begin{minipage}[b]{75mm}
\centerline{\includegraphics{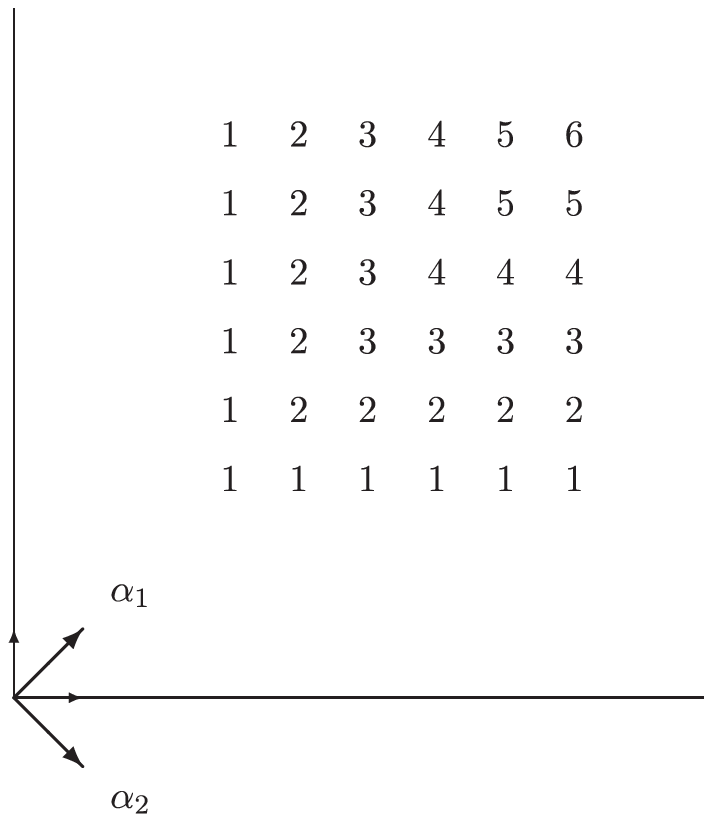}}
\caption{}\label{fig2}
\end{minipage}
\end{figure}

$ K^{H_2}= K \cap L$ is isomorphic to $U(2)$. We assume that the
positive root of $(\bt, \bk^{H_2})$ is $\alpha _2$ whereas the
$\alpha_1$ can be identif\/ied with a character of the center of
$K^{H_2}$. So the restriction of each $K$-type with highest
weight $(m+3)\alpha _1$ to $ K^{H_2}$ is a sum of characters
$d$$\alpha_1$ with $ -(m+3)\leq d \leq (m+3)$. Hence $A_\bq $ is
not $K^{H_2} $-f\/inite and thus by T.~Kobayashi's theorem $A_\bq $
is not a discrete sum of irreducible representations of $H_2$.
This proves the second claim of Proposition~\ref{propositionI.1}.

 The groups $H_1$ and $H_2$ are conjugate under the outer
automorphism which changes the sign of $e_2$. Hence the simple
positive root of $(\bt, \bk^{H_1})$
 can be
identif\/ied with $\alpha _1$ and $\alpha_2$ can be identif\/ied with
a character of the center of $K^{H_1}$ (see Fig.~\ref{fig1}). A $K$-type with highest
weight $\Lambda +m_1(e_1 +e_2)+ 2m_2 \ e_1$ is a tensor product of
a representation with highest weight $(3+m_1+m_2,3+m_1+m_2)$ and a
representation with highest weight $(m_2,-m_2)$. Its restriction
to $K ^{ H_1}$ is a direct sum of representations with highest
weights $(3+m_1+m_2 +i,3+m_1+m_2-i)$, $-m_2 \leq i \leq m_2$.
Fig.~\ref{fig2} shows  the highest weights of the $K^{ H_1}$-types for
the restriction of $A_\bq $ to $K ^{H_1}$. Their multiplicities
are indicated by a number.

%\bigskip

%\setlength{\unitlength}{0.7cm}
%\begin{picture}(10,15)
% \thicklines
% \put(1.0,4.0){\vector(1,1){1.0}}
% \put(2.4,5.4){$\alpha_1$}
% \put(1.0,4.0){\vector(1,-1){1.0}}
% \put(2.4,2.4){$\alpha_2$}
% \thinlines
% \put(1.0,4.0){\vector(0,1){1.0}}
% \put(1.0,4.0){\vector(1,0){1.0}}
% \put(1.0,4.0){\line(1,0){10.0}}
%\put(1.0,4.0){\line(0,1){10.0}}
%
% \multiput(4.0,7.0)(1.0,0){6}{1}
% \multiput(5.0,8.0)(1.0,0){5}{2}
% \multiput(6.0,9.0)(1.0,0){4}{3}
% \multiput(7.0,10.0)(1.0,0){3}{4}
% \multiput(8.0,11.0)(1.0,0){2}{5}
% \multiput(9.0,12.0)(1.0,1.0){1}{6}
% \multiput(4.0,8.0)(0,1.0){5}{1}
% \multiput(5.0,9.0)(0,1.0){4}{2}
% \multiput(6.0,10.0)(0,1.0){3}{3}
% \multiput(7.0,11.0)(0,1.0){2}{4}
% \multiput(8.0,12.0)(0,1.0){1}{5}
% \put(5,1){Figure 2}
%\end{picture}

Thus $A_\bq$ is $K^{H_1}$-f\/inite. This completes the proof of the
f\/irst assertion of Proposition~\ref{propositionI.1}.

\medskip \noindent {\bf Remark.}   A second series of
representations is obtained if we def\/ine another $\theta $-stable
parabolic subalgebra $\bq'$ using the matrix $Q_1\in \bg$ instead
of $Q_2$. We obtain a representation $A_{\bq'}$ which is not
equivalent to $A_{\bq}$. The same arguments as in the previous
case prove that restriction of $A_{\bq'}$ to $H_1$ does not have a
purely discrete spectrum whereas the restriction to $H_2$ is a
direct sum of irreducible unitary representations.

The representation ${\mathcal A}_\bq$ of $SL^{+/-}(4,\bR)$
(determinant $\pm 1$) obtained by inducing representation~$A_\bq$  is irreducible and
its restriction to $SL(4,\bR)$ is equal to $A_\bq \oplus
A_{\bq'}$. Hence the restriction of ${\mathcal A}_\bq $ to $H_1$
does not have discrete spectrum.

\section[The restriction of $A_{\bq}$ to the symplectic group $H_1$]{The restriction of $\boldsymbol{A_{\bq}}$ to the symplectic group $\boldsymbol{H_1}$}\label{sec3}

In this section we determine $H_1$-types  of $A_\bq$. Our
techniques are based on homological al\-geb\-ra and the construction
of an ``enlarged complex'' whose cohomology computes the
restriction. We introduce it  in {\bf 4.1} %III.1
for semisimple connected
Lie groups~$H$ and connected reductive subgroups~$H$. Then we will
compute the restriction to~$H_1$ by restricting~$A_\bq$ to a
subgroup conjugate to~$H_1$. The motivation for this ``enlarged
complex'' or ``branching complex'' is the same as when one is
restricting holomorphic functions to a complex submanifold, and
identifying the functions with their normal derivatives along the
submanifold. In our case we are working with (formalizations of)
dif\/ferential forms satisfying a similar dif\/ferential equation, so
it is natural to try to identify them with their ``normal
derivatives''; this is what is formalized in our def\/inition. As it
turns out, with the appropriate conditions (well aligned parabolic
subgroups, vanishing of the cohomology in many degrees, and the
non-vanishing of explicit classes corresponding to small
$K$-types) we can indeed make the calculation of the branching law
ef\/fective, at least in the examples at hand.

\medskip

\noindent{\bf 4.1.} We def\/ine  another complex for a semisimple
connected Lie group $G$ and a connected reductive subgroup $H$
satisfying the assumptions of {\bf 2.3}.
Let $\bC_{\lambda_H}$ be the
one dimensional representation of $L^H$ def\/ined by $\lambda
\otimes \wedge ^{{\rm top}}(\bu \cap \bs)$. Then
\[\bC_\lambda^\sharp = \bC_\lambda \otimes \wedge^{{\rm top}}\bu
= \bC_{\lambda_H}\otimes \wedge ^{{\rm top}}\bu^H
=\bC_{\lambda_H}^\sharp.\]
Consider the complex $L_H^*$
\begin{equation}\label{bigcomplex}(\mbox{Hom}_{L\cap K\cap H}(
U(\bg),\mbox{Hom}(\wedge^i r_H ,\bC_{\lambda_H} ^\sharp ))_{K\cap
H},d_H).\end{equation} Here $r_H=r_G \cap \bh$ and $d_H$ is
def\/ined analogously to the dif\/ferential $d$ in {\bf 2.2}. %I.2.
 As a left $U(\bl^H)$-module
\[ U(\bg) = Q \otimes U(\bh ),
\]
where $Q $ is the symmetric algebra $S( \bs)$. (See \cite[2.56]{K-V}.) We have
\begin{gather*}
\mbox{Hom}_{L\cap K\cap H}( U(\bg),\mbox{Hom}(\wedge^i
r_H ,\bC_{\lambda_H} ^\sharp ))_{K\cap H}
 =
 \mbox{Hom}_{L\cap K\cap H}(
Q\otimes U(\bh ), \mbox{Hom}(\wedge^i r_H,\bC_{\lambda_H} ^\sharp ))_{K\cap H} \\
\qquad{} =   \mbox{Hom}_{L\cap K\cap H}( U(\bh),
\mbox{Hom}(\wedge^i r_H ,Q^* \otimes\bC_{\lambda_H} ^\sharp
))_{K\cap H}.
\end{gather*}
$U(\bg)$ acts on the enlarged complex from the right and a quick
check shows that $d_H$ also commutes with this action and
therefore we have an action of $U(\bg)$ on the cohomology of the
complex.

We also consider the ``large  complex'' $L^*$
\[
(\mbox{Hom}_{K \cap L \cap
H}( U(\bg),\mbox{Hom}(\wedge^i r_G, \bC_\lambda ^\sharp ))_{K\cap
H},\,d).
 \]

We have $r_G = r_H \oplus (\bu \cap \bs) \oplus ( p_L  \cap \bs)$,
and so
\[\wedge^i r_G = \oplus _{l+k=i} \wedge ^k r_H \otimes \wedge ^l
 (\bu \cap \bs \oplus p_L \cap \bs) .\]

Using this we can def\/ine a ``pull back''  map of forms
\begin{gather*}  \mbox{\bf pb}^i_H: \ \
\mbox{Hom}_{L\cap K\cap H}( U(\bh),\mbox{Hom}(\wedge^i
r_H,Q^*\otimes \bC_{\lambda_H} ^\sharp ))_{K \cap H}
\\  \phantom{\mbox{\bf pb}^i_H: \ \ }  \rightarrow
\mbox{Hom}_{K \cap L \cap H}( U(\bg),\mbox{Hom}(\wedge^i r_G,
\bC_\lambda ^\sharp ))_{K\cap H}.
\end{gather*}
The pullback map commutes with the right action of $U(\bh )$ and
induces a map of complexes.
This is the main observation we use in order to analyze the action
of $H$ on the cohomologically induced module.

\medskip
\noindent {\bf 4.2.}  For the rest of this section we assume
that  $G=SL(4,\bR)$. We will show that there exists a symplectic
subgroup, which we denote by $H_1^w$
 conjugate to $H_1$ by an element $w$, so that
 the pullback induces a nontrivial map in cohomology.
  Since the restriction of $A_\bq $ depends only on the conjugacy
class of $H_1$ this determines the restriction.

 Since
\[
\left(\begin{array}{cccc} 1&0&0&0 \\ 0&0&-1&0 \\0&1&0&0 \\0&0&0&1
\end{array}\right)Q_1
\left(\begin{array}{cccc} 1&0&0&0 \\ 0&0&1&0 \\0&-1&0&0 \\0&0&0&1
\end{array}\right) = \left(\begin{array}{cccc} 0&0&-1&0 \\ 0&0&0&-1\\1&0&0&0
\\0&1&0&0
\end{array}\right),
\]
by abuse of notation we will also write $H_1$ and $H_1'$ for the
groups def\/ined by the skew symmetric form
\[\left(\begin{array}{cccc} 0&0&-1&0 \\ 0&0&0&-1\\1&0&0&0
\\0&1&0&0
\end{array}\right).\]
Thus{\samepage  \[\bh_1= \left(
\begin{array}{cc} A& X\\ Y & -A^{\rm tr}
\end{array}\right)\]
for symmetric matrices $X$ and $Y$.}

\medskip Recall that $\bq $ is def\/ined by
\[i \ Q_2=\left(
\begin{array}{cc}i\  J & 0 \\ 0 & i\ J \end{array}\right) \in \bg\]
and that
 $\bg = \bh_1 \oplus \bs_1$.
We need the f\/ine structure of the parabolic relative to the
symmetric subgroup, in order to compare the cohomology of these
complexes during the branching.

\begin{mylemma}\label{lemma4.1}
Under the above assumptions

\begin{enumerate}
\itemsep=0pt
\item[a)] $\bl \cap \bh_1 $ is isomorphic to $sl(2,\bR )
\oplus i\bR $ and   $\dim (\bu \cap \bh_1) =3$;

\item[b)] the representation of $L\cap H_1$ acts by a
nontrivial character $\mu_{1}$ with differential $(e_1+e_2)$ on
the one dimensional space $\bu \cap \bs_1 $;

\item[c)]  $\bl \cap \bs_1$ is a direct sum of the
trivial representation and the adjoint representation of $\bl \cap
\bh_1 $.

\item[d)]  $\bu \cap \bk= \bu \cap \bk \cap \bh_1$ has
dimension $1$.
\end{enumerate}
\end{mylemma}

\begin{proof} We have
\[\blo \cap \bh_1
=\left(\begin{array}{cccc} a&b&x&0
\\ -b&a&0&x\\y&0&-a&b
\\0&y&-b&-a
\end{array}\right).\] The nilradical of a parabolic subalgebra with this
Levi subalgebra has dimension 3.

The dimension of $\bl \cap \bh_1 \cap \bk$ is 2. Hence the
dimension of $\bu \cap \bk \cap \bh_1$ is 1. On the other hand the
dimension of $\bl \cap \bk $ is 4. So the dimension of $\bu \cap
\bk$ is 1. Since $\bu \cap \bk \cap \bh_1 \subset \bu \cap \bk$ we
have equality.

$\bu \cap \bs_1 $ is in the roots spaces for roots $e_1+e_2 +(f_1
-f_2)$ and $e_1+e_2-(f_1-f_2)$.  Hence $\bl \cap \bh_1 \cap \bk$
acts on $\bu \cap \bs_1 $ by $e_1 +e_2$.

$\bl \cap \bh_1 $ acts on the 4 dimensional space $\bl \cap \bs_1$
via the adjoint representation.
\end{proof}

The representation of $L\cap H_1$ on the symmetric algebra
$ S((\bu \cap \bs_1) \oplus (\overline{\bu} \cap \bs_1) \oplus
(\bl \cap \bs_1))$ is isomorphic to a direct sum of
representations $\mu_{1}^{n_1}\otimes \mu_{1}^{-m_1}\otimes
\mbox{ad} ^{r_1 }$ with $n_1, m_1 , r_1 \in \bN$.
These powers of $\mu_1$ will label the constituents in the
branching law; it will also be sometimes  convenient to think of
their dif\/ferentials in additive notation.

Now it is important to note, that the parameter
$\lambda_{H_1}\otimes \mu_{1}^{n_1}$,  $0 \leq n_1 $ is in the
good range~\cite{K-V} and thus the representation on the
cohomology in degree $1 = \dim (\bu \cap \bk \cap \bh_1)$
of the  complex~$L_{H_1}^*$ has composition factors isomorphic to
\[
A_{\bq \cap {\bh_1}}(\lambda_{H_1}\otimes \mu_{1}^{n_1}),
\] where $0 \leq n_1 $.
In particular $A_{\bq \cap {\bh_1} }(\lambda_{H_1} )$ is an
$(\bh_1 , K\cap H_1)$-submodule module of the cohomology of
$L_{H_1}^*$.

\begin{myprop}\label{proposition4.2}
 $A_{\bq \cap \bh_1}(\lambda_{H_1})$ is a composition factor
 of the restriction of $A_\bq $ to $(\bh_1,K\cap H_1)$.
\end{myprop}

\begin{proof} Note that $\dim (\bu \cap \bk)= \dim (\bu
\cap \bk \cap \bh_1)$ and that 1 = $\dim (\bu \cap \bk)$ is the
degree  in which the complexes def\/ining the representations
$A_\bq$ respectively $A_{\bq \cap \bh_1}(\lambda_{H_1})$ both have
nontrivial cohomology \cite{K-V}. Considering the complex def\/ining
$A_\bq$ as a subcomplex of the ``large complex''~$L^*$ the pullback
$\mbox{\bf pb}_{H_1}^i $ of forms def\/ines a $(\bh_1,K\cap
H_1)$-equivariant map \[A_\bq \rightarrow \oplus _{n_1=0}^\infty
A_{\bq \cap \bh_1}(\mu_1^{n_1}\otimes \lambda_{H_1}).\]

Recall the def\/inition of the $K$-module ${\mathcal
R}_K^s(\lambda)$ from V.5.70 in \cite{K-V}. We have bottom layer
maps of $\bko $-modules.
\[{\mathcal B}(\lambda): A_\bq \rightarrow {\mathcal R}_K^1(\lambda_0)  \]
and
\[{\mathcal B}(\lambda_{H_1}):A_{q \cap
\bh_1}(\lambda_{H_1} ) \rightarrow {\mathcal R}_{K \cap
H_1}^1(\lambda_{H_1}),
 \] where $\lambda_0$ is the trivial
character of $L\cap K$. These maps are def\/ined by the inclusion of
of complexes and hence of forms.  See Theorem~V.5.80 and its proof
in~\cite{K-V}. The minimal $K$-types of $A_\bq$, respectively
$K^{H_1}$-type of $ A_{\bq \cap \bh_1}(\lambda_{H_1})$ are in the
 bottom layer.

On the other hand we have an inclusion of complexes
(the notation in analogy with the case in~{\bf 4.1}, now for
the case where we take $G = K$)
\begin{gather*} \mbox{\bf pb}_{H_1\cap K}^i : \ \ \mbox{Hom}_{K \cap L \cap H_1}(
U(\bk),\mbox{Hom}(\wedge^i (r_G\cap \bk \cap \bh_1), \bC_{\lambda
_{H_1 }} ^\sharp
))_{K \cap H_1} \\
\phantom{\mbox{\bf pb}_{H_1\cap K}^i :} \ \ \rightarrow   \mbox{Hom}_{K \cap L \cap H_1}(
U(\bk),\mbox{Hom}(\wedge^i (r_G\cap \bk), \bC_{\lambda
_{0}}^\sharp ))_{K \cap H_1}  .
\end{gather*}
But $ \bC_{\lambda _{0}}^\sharp = \bC_{\lambda _{H_1}} ^\sharp $
and $r_G\cap \bk \cap \bh_1 = r_G\cap \bk $ and so using a forgetful
functor  we may consider
\begin{gather*}
\mbox{Hom}_{K \cap L }(
U(\bk),\mbox{Hom}(\wedge^i (r_G\cap \bk), \bC_{\lambda
_{0}}^\sharp ))_{K }
\end{gather*}
 as a subspace, respectively subcomplex, of
\begin{gather*}
 \mbox{Hom}_{K \cap L \cap H_1}(
U(\bk),\mbox{Hom}(\wedge^i (r_G\cap \bk \cap \bh_1), \bC_{\lambda
_{0}}^\sharp ))_{K \cap
H_1} \\
\qquad{}=\mbox{Hom}_{K \cap L \cap H_1}( U(\bk\cap
\bh_1),\mbox{Hom}(\wedge^i (r_G\cap \bk \cap \bh_1),Q_H \otimes
\bC_{\lambda _{0}}^\sharp ))_{K \cap H_1},
\end{gather*}
where $Q_H$ is the symmetric algebra of the complement of $\bh_1
\cap \bk $ in $\bk$. Note that $K\cap H_1$ and $K$~is again a
symmetric a pair and so we have a bottom-layer map for the
representation
\[{\mathcal
R}_K^1(\lambda_0) \rightarrow {\mathcal R}_{K \cap
H}^1(\lambda_{H_1}).\] Since the representation ${\mathcal
R}_K^1(\lambda_0)$ is irreducible restricted to $K\cap H_1$ this
map is an isomor\-phism.~~
\end{proof}

\noindent {\bf Def\/inition.}  We call
 $A_{\bq \cap \bh_1}(\lambda_{H_1})$ the minimal $H_1$-type of
 $A_\bq$.

\begin{mytheorem}\label{theorem4.3}
The representation $ A_\bq$ restricted to $H_1$ is the direct sum
of the representations each occuring with multiplicity one, namely
\[{A_\bq}_{|H_1}= \oplus _{n_1=0}^\infty A_{\bq
\cap \bh_1}(\mu_1^{n_1}\otimes \lambda_{H_1}).\]
\end{mytheorem}

\begin{proof}  By the proof of the lemma $A_{\bq \cap
\bh_1}(\lambda_{H_1} )$ is a submodule of the restriction of
$A_\bq$ to the symplectic group $H_1$. Its minimal $K^{H_1}$-type
is also a minimal $K$-type of $A_\bq$ and hence occurs with
multiplicity one. Hence $A_{\bq \cap \bh_1}(\lambda_{H_1} )$ is a
$H_1$-type of $A_\bq$ with multiplicity one.

  The minimal $K^{H_1}$-type of $A_{\bq \cap \bh_1}(\lambda)$
has highest weight $\lambda+3e_1 +3e_2$. The roots of $\bu \cap
\bh_1 \cap p$ are $2e_1$, $2e_2$. Applying successively the root
vectors to the highest weight vector of the minimal
$K^{H_1}$-type of $A_{\bq \cap \bh_1}(\lambda)$ we deduce that
 $A_{\bq \cap \bh_1}(\lambda )$ contains  the $K^{ H_1}$-types with
highest weight $((3+2r_1)e_1 +(3+2r_2)e_1 +\lambda)$, $r_1,r_2 \in
\bN$. Theorem 8.29 in \cite{K-V} show that all these
$K^{H_1}$-types have multiplicity one.  Fig.~\ref{fig3} shows the
$K^{H_1}$-type multiplicities of $A_\bq(\lambda_{H_1})$.

Note that we are here using quite a bit of a priori information
about the derived functor modules for the smaller group; on the
other hand, the branching problem has essentially been reduced to
one for compact groups, $K$-type by $K$-type.

\begin{figure}[t]
\centerline{\includegraphics{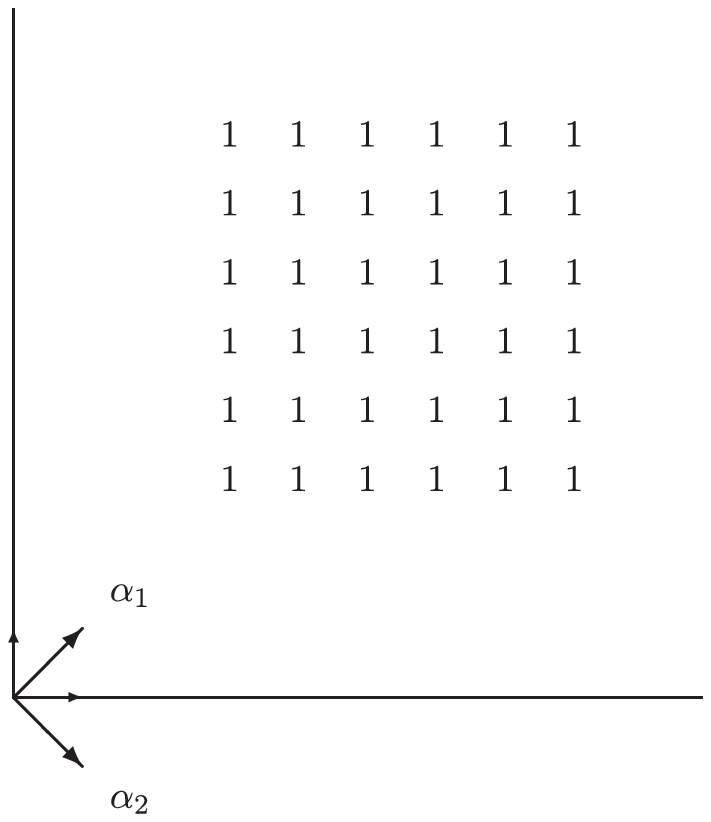}}

\caption{}
\label{fig3}
\end{figure}

The Borel subalgebra of $\bk \cap \bh _1 $ acts on the one
dimensional space $\bu \cap \bs_1 $ by a character $\mu_1$ with
dif\/ferential $(e_1+e_2)$. Let $Y \not= 0$ be in $\bu \cap \bs_1 $
and $v\not =0$ a highest weight vector of the minimal $K$-types of
$A_\bq.$ Then $Y^n\cdot v \not= 0$ is also the highest weight of
an $K^{H_1}$-type of highest weight $(3+n)e_1+(3+n)e_2$ of
$A_\bq$.

Let $X_k \not= 0$ be  in $\bu\cap \bk$. The linear map
\[T_s:
U(\bg ) \rightarrow \wedge^1 r_G \otimes \bC_{\lambda_0}^\sharp\]
which maps 1 to $X_k \otimes \bC_{\lambda_0}^\sharp$ is non-zero
in cohomology and its class $[T_s]$ is the highest weight vector
of the minimal $K$-type. But
\begin{eqnarray*} Y \cdot T_s &\in & \mbox{Hom}_{L\cap K\cap H_1}(
\bs_1\otimes U(\bh _1 ), \mbox{Hom}(\wedge^s
r_{H_1},\bC_{\lambda_{H_1}} ^\sharp )),
\end{eqnarray*}
so may consider
\begin{eqnarray*}
 Y \cdot T_s &\in &\mbox{Hom}_{L\cap K\cap
H_1}( U(\bh_1), \mbox{Hom}(\wedge^i r_{H_1} ,\bs_1^*
\otimes\bC_{\lambda_{H_1}} ^\sharp ))_{K\cap H_1}.
\end{eqnarray*}
Hence  $0 \not = [Y \cdot T_s] = Y \cdot [T_s] \in A_{\bq \cap
\bh_1}(\mu_1 \otimes\lambda_{H_1}) $
 and thus $A_{\bq \cap
\bh_1}(\lambda_{H_1} \otimes \mu_1)$ is an $H_1$-type  of
$A_\bq.$ The same argument shows that $A_{\bq \cap
\bh_1}(\lambda_{H_1} \otimes\mu_1^n), \ n\in \bN,$ is a $H_1$-type
of $A_\bq$.

Now every $K$-type with highest weight $(n,n)$ has multiplicity
$n-2$ and is contained in exactly $n-2$ composition factors. The
multiplicity computations in Section~\ref{sec2} now show that every
composition factor is equal to $A_{\bq \cap \bh_1}(\lambda_{H_1}
\otimes \mu_1^n)$ for some $n$. See Fig.~\ref{fig4}.
\end{proof}

 \noindent {\bf Remark.} Another proof of Theorem~\ref{theorem4.3} can be
 obtained using Proposition \ref{proposition5.1} %IV.1
 and the ideas of {\bf 6.1}.

 \begin{figure}[t]
\centerline{\includegraphics{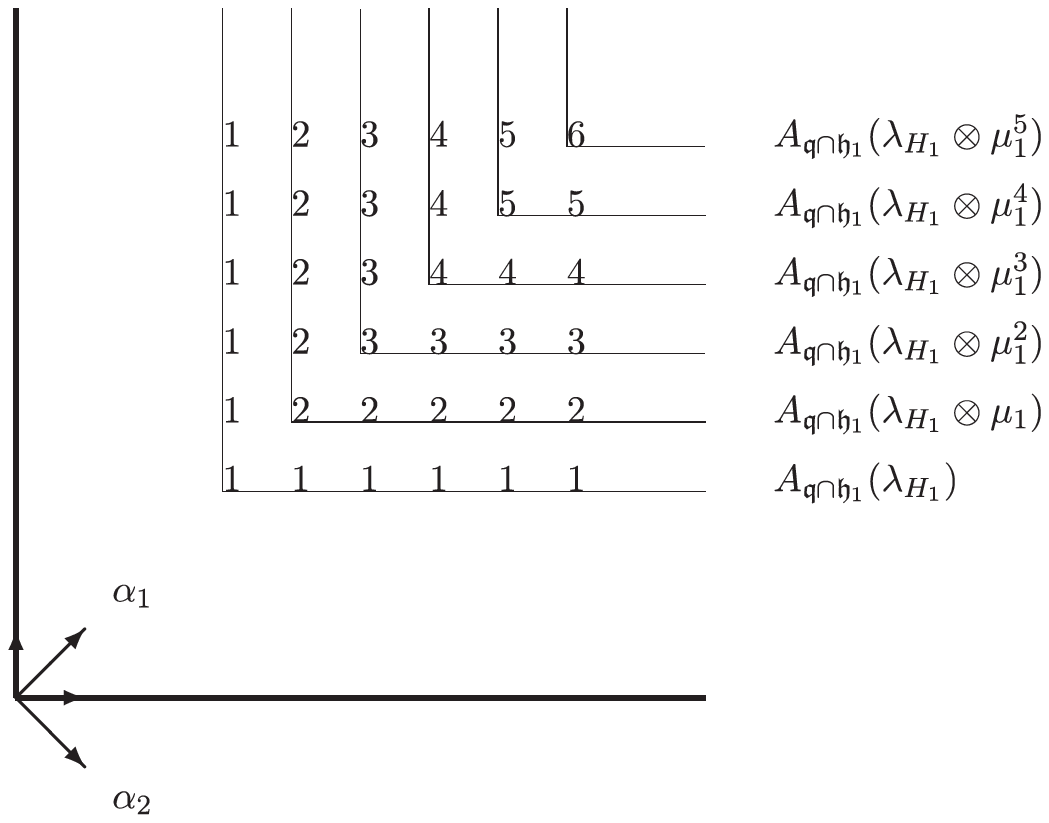}}

\caption{}
\label{fig4}
\end{figure}

%\bigskip
%\setlength{\unitlength}{0.7cm}
%\begin{picture}(10,15)
% \thicklines
% \put(1.0,4.0){\vector(1,1){1.0}}
% \put(2.4,5.4){$\alpha_1$}
% \put(1.0,4.0){\vector(1,-1){1.0}}
% \put(2.4,2.4){$\alpha_2$}
% \linethickness{0.5mm}
% \put(1.0,4.0){\vector(0,1){1.0}}
% \put(1.0,4.0){\vector(1,0){1.0}}
%
% \put(1.0,4.0){\line(1,0){10.0}}
% \put(1.0,4.0){\line(0,1){10.0}}
%
% \multiput(4.0,7.0)(1.0,0){6}{1}
% \multiput(5.0,8.0)(1.0,0){5}{2}
% \multiput(6.0,9.0)(1.0,0){4}{3}
% \multiput(7.0,10.0)(1.0,0){3}{4}
% \multiput(8.0,11.0)(1.0,0){2}{5}
% \multiput(9.0,12.0)(1.0,1.0){1}{6}
% \multiput(4.0,8.0)(0,1.0){5}{1}
% \multiput(5.0,9.0)(0,1.0){4}{2}
% \multiput(6.0,10.0)(0,1.0){3}{3}
% \multiput(7.0,11.0)(0,1.0){2}{4}
% \multiput(8.0,12.0)(0,1.0){1}{5}
%
% \linethickness{0.03mm}
%
% \put(4.0,7.0){\line(1,0){7.0}}
%  \put(4.0,7.0){\line(0,1){7.0}}
%  \put(5.0,8.0){\line(1,0){6.0}}
%  \put(5.0,8.0){\line(0,1){6.0}}
%  \put(6.0,9.0){\line(1,0){5.0}}
%  \put(6.0,9.0){\line(0,1){5.0}}
%  \put(7.0,10.0){\line(1,0){4.0}}
%  \put(7.0,10.0){\line(0,1){4.0}}
%  \put(8.0,11.0){\line(1,0){3.0}}
%  \put(8.0,11.0){\line(0,1){3.0}}
%  \put(9.0,12.0){\line(1,0){2.0}}
%  \put(9.0,12.0){\line(0,1){2.0}}
%
%  \put(12,7){$A_{\bq \cap \bh_1}(\lambda_{H_1})$}
%  \put(12,8){$A_{\bq \cap \bh_1}(\lambda_{H_1}\otimes\mu_1)$}
%  \put(12,9){$A_{\bq \cap \bh_1}(\lambda_{H_1}\otimes \mu_1^2)$}
%  \put(12,10){$A_{\bq \cap \bh_1}(\lambda_{H_1}\otimes\mu_1^3)$}
%  \put(12,11){$A_{\bq \cap \bh_1}(\lambda_{H_1}\otimes\mu_1^4)$}
%  \put(12,12){$A_{\bq \cap \bh_1}(\lambda_{H_1}\otimes\mu_1^5)$}
% \put(5,1){Figure 4}
%\end{picture}

By Proposition~8.11 in \cite{K-V} for any $(\bg,K)$-module X we have
 \[
 \mbox{Hom}_{\bg,K}(X,A_\bq) = \mbox{Hom}_{\bl,K \cap
 L}(H_s(\bu,X),\bC^\sharp),
 \]
 where $H_s(\bu,X)$ is the Lie algebra homology as def\/ined in \cite{K-V} and
 $s$ = dim($\bu \cap \bk$).
Thus we have a ``Blattner type formula'' for the $H_1$-types of
$A_\bq$.

\begin{mycor}\label{corollary4.4}
Let $V$ be an irreducible $(\bh_1, K^{H_1})$-module. Then
\begin{gather*}
{\rm dim \,Hom}_{\bh_1,K\cap H_1} ( V,A_\bq ) =
  \sum_i {\rm dim \, Hom}_{(\bl \cap \bh_1,K^{H_1}\cap L)}
(H_1(\bu \cap \bh_1,V), S^i(\bu \cap \bs_1 )\otimes
\bC_{H_1}^\sharp).
\end{gather*}
\end{mycor}

 \noindent {\bf Remark.}
The $H_1(\cdot,\cdot)$ on the right refers to homology in degree one.
The maximal Abelian split subalgebra $\ba_1$  in $\bl \cap \bh_1$
are the diagonal matrices. So the parabolic subgroup of the
Langlands parameter of the $H_1$-types of $A_\bq$ is the so
called ``mirabolic'', i.e.\ the maximal parabolic subgroup with
Abelian nilradical. The other  parts of Langlands parameter can de
determined using the algorithm in~\cite{K-V}.

We consider in Theorem~\ref{theorem4.3} the restriction of a ``small'' representation
$A_\bq$ of Spin(3,3) to Spin(2,3)  similar to the restriction of
``small'' discrete series representations of $SO(n+1,n)$ to $SO(n,n)$
considered by B.~Gross and N.~Wallach in~\cite{varfa}. It would be
interesting to see if their techniques could be adapted to the
problem discussed in the paper.

\section[Restriction of $A_\bq$ to the group $H_1'$]{Restriction of $\boldsymbol{A_\bq}$ to the group $\boldsymbol{H_1'}$}\label{sec4}

In this section we describe $H_1'$-types of $A_\bq$ using the
same techniques as in the previous section.

\medskip \noindent
{\bf 5.1.}  For $H_1'$  we consider the  complex
\begin{equation*}\label{bigcomplex}(\mbox{Hom}_{L\cap K\cap H_1'}(
U(\bg),\mbox{Hom}(\wedge^i r_{H_1'} ,\bC_{\lambda_{H_1'}} ^\sharp
))_{K\cap H_1'},\, d_{H_1'}) \end{equation*}
 and the map
\begin{gather*}  \mbox{\bf pb}_{H_1'}^i: \ \
\mbox{Hom}_{L\cap K\cap H_1'}( U(\bh _1'),\mbox{Hom}(\wedge^i
r_{H_1'},Q^*\otimes \bC_{\lambda_H} ^\sharp )_{K\cap H_1'})_{K
\cap H_1'}
\\\phantom{ \mbox{\bf pb}_{H_1'}^i: \ \ }{} \rightarrow
\mbox{Hom}_{K \cap L \cap H_1'}( U(\bg),\mbox{Hom}(\wedge^i r_G,
\bC_\lambda ^\sharp ))_{K \cap H_1'}.
\end{gather*}

 We write $\bg = \bh_1' \oplus \bs_1'$. The intersection
$\bu \cap \bs_1'$ is 2-dimensional and the representation of the
group of $L\cap H_1'$ on $\bu \cap \bs_1' $ is reducible and thus
a sum of 2 one dimensional representations $\chi_1 \oplus \chi_2$.
The weights of these characters  are $2e_1$ and $2e_2$. So the
symmetric algebra $S(\bu \cap \bs_1' )$ is a direct sum of one
dimensional representations of $L \cap H_1' $ with weights $ 2m_1
e_1 +2m_2 e_2$.

In the cohomology in degree 1 of the  complex $L_{H_1'}^*$ we have
composition factors
\[A_{\bq\cap \bh_1'}(\lambda_{H_1'}\otimes \chi_1^{n_1}\otimes \chi_2^{n_2})\]
with $0 \leq n_1, n_2 $. In particular
 $A_{\bq \cap \bh_1'}(\lambda_{H_1'} )$ is an $(\bh_1', K\cap H_1')$-submodule
module of the cohomo\-logy in degree 1.

\begin{myprop}\label{proposition5.1}
$A_{\bq \cap \bh_1'}(\lambda_{H_1'} )$ is a composition factor of
the restriction of $A_\bq$ to $(\bh_1',H_1' \cap K)$.
\end{myprop}

\begin{proof} The maximal compact subgroups of $H_1$ and
$H_1'$ are identical. Thus $\dim  \bu \cap \bk \cap \bh_1=
\dim \bu \cap \bk \cap \bh_1=1$ and the minimal $K$-type is
irreducible under restriction to $K^{H_1'}$. Thus the same
argument as in Lemma~\ref{lemma4.1} completes the proof.
\end{proof}

\noindent {\bf Def\/inition.} We call $A_{\bq \cap
\bh_1'}(\lambda_{H_1'} )$ the minimal $H_1'$ type of $A_\bq$.

\begin{mytheorem}
The representation $ A_\bq$ restricted to $H_1'$ is the direct sum
of the representations each occurring with multiplicity one,
namely
\[{A_\bq}_{|H_1'}= \oplus _{n_1,n_2=0}^\infty
A_{\bq\cap \bh_1'}(\lambda_{H_1'} \otimes \chi_1^{n_1}\otimes
\chi_2^{n_2}).\] Their minimal $K\cap H_1'$ -types have highest
weights $(3+m_1+m_2 +i,3+m_1+m_2-i)$, $-m_2 \leq i \leq m_2$.

\end{mytheorem}

\begin{proof} The proof is the same in the previous
section where we proved that the representations $A_{\bq\cap
\bh_1'}(\lambda_{H_1'} \otimes \chi_1^{n_1}\otimes \chi_2^{n_2})$,
$n_1, n_2 \in \bN$ appear in the restriction of the $A_\bq $ to
$H_1'$

The $K\cap H_1'$ -types of all unitary representations of
$GL(2,\bC )$  have multiplicity one. If the  minimal $K\cap H_1'$-type has highest weight $l_1e_1+ l_2e_2-2$, then the highest
weights of the other $K\cap H_1'$-types are $(l_1+j)e_1+
(l_2+j)e_2$. Multiplicity considerations of $K\cap H_1'$-types of
$A_\bq $ conclude the proof.
\end{proof}

 \noindent
 {\bf Remark.} Using Proposition \ref{proposition4.2}, {\bf 4.1} and the ideas of {\bf 6.1}
we can obtain another proof of this theorem.

\medskip

Fig.~\ref{fig5} shows the decomposition into irreducible representations.
The highest weights of the $K\cap H_1'$-types of a composition
factors lie on the lines. For each highest weight there is exactly
one composition factor which has a $K\cap H_1'$-type with this
weight as a minimal $K\cap H_1'$-type.

We have again a ``Blattner-type formula'' for the $H_1'$ types.

\begin{mycor}
Let V be an irreducible $(\bh_1',(K \cap H_1'))$-module, then
\begin{gather*} \dim {\rm Hom}_{\bh_1',K\cap H_1'} (
V,A_\bq)=
\sum _i \dim {\rm Hom}_{(\bl \cap \bh_1'),K\cap H_1'
\cap L} (H_1(\bu \cap \bh_1',V), S^i(\bu \cap \bs_1' )\otimes
\bC_{H_1'}^\sharp).
\end{gather*}
\end{mycor}

\noindent {\bf 5.2.}  All the $H_1'$-types $A_{\bq \cap
\bh_1'}(\lambda_{H_1'}\otimes \chi_1^{n_1}\otimes \chi_2^{n_2})$
of $A_\bq$ are simply unitarily induced principal series
representations of $GL(2,\bC)$.

\section{Pseudo dual pairs}\label{sec5}

\noindent {\bf 6.1.}  Suppose now that $G$ is a reductive
connected Lie group with maximal compact sub\-group~$K$, Cartan
involution $\theta$ of $\bg$. For an involution
 \[\tau: G\rightarrow G\]
 commuting with
$\theta$ we def\/ine
\[\tau' = \tau \circ \theta\qquad \mbox{and}\qquad H = G^{\tau},\qquad H' = G^{\tau'}.\]

\noindent {\bf Def\/inition.} We call $H$ and $H'$ a {\it pseudo dual
pair}.

\medskip Since $\theta$, $\tau$ and $\tau'$ commute

\begin{mylemma}
Suppose that $H$ and $H'$ are a pseudo dual pair in $G$. Then
\begin{enumerate}\itemsep=0pt
\item[\rm 1)]  $K \cap H = K\cap H'$;
 \item[\rm 2)]  we have $\bp = \bp^{\tau} \oplus
\bp^{\tau'}$;
 \item[\rm 3)]
$U(\bg) \cong S(\bp^{\tau'}) S(\bp^{\tau}) U(\bk)$ and so we can
write every element in $U(\bg)$ as a linear combination of terms
of the form $\omega_2 \omega_1 \omega_k $ where $\omega_k \in
U(\bk)$, $\omega_1 \in U(\bh)$ and $\omega_2 \in U(\bh')$.
\end{enumerate}
\end{mylemma}

Suppose now that $\pi $ be an irreducible unitary
$(\bg,K)$-module.  Then $\pi $ is $K \cap H$-admissible if and
only if it is $K \cap H'$-admissible and by Theorem 4.2 of~\cite{ko1}
its restriction to $\bh$ is a direct sum  of irreducible
$(\bh,K\cap H)$-modules  if and only if $\pi $ is a direct sum of
irreducible $(\bh',K\cap H')$-modules.

Suppose that  $\pi $ is  $H\cap K$ admissible  and that $V_k
\subset \pi $ is a minimal $K$-type of $\pi$. If $\pi_o$ is a~$(\bh,K\cap H)$-module, which occurs in the restriction of $\pi$
to $\bh$ then
\[ \pi_o \cap U(\bh')V_k \neq 0. \]

\begin{figure}[t]
\centering
\begin{minipage}[b]{75mm}
\centerline{\includegraphics{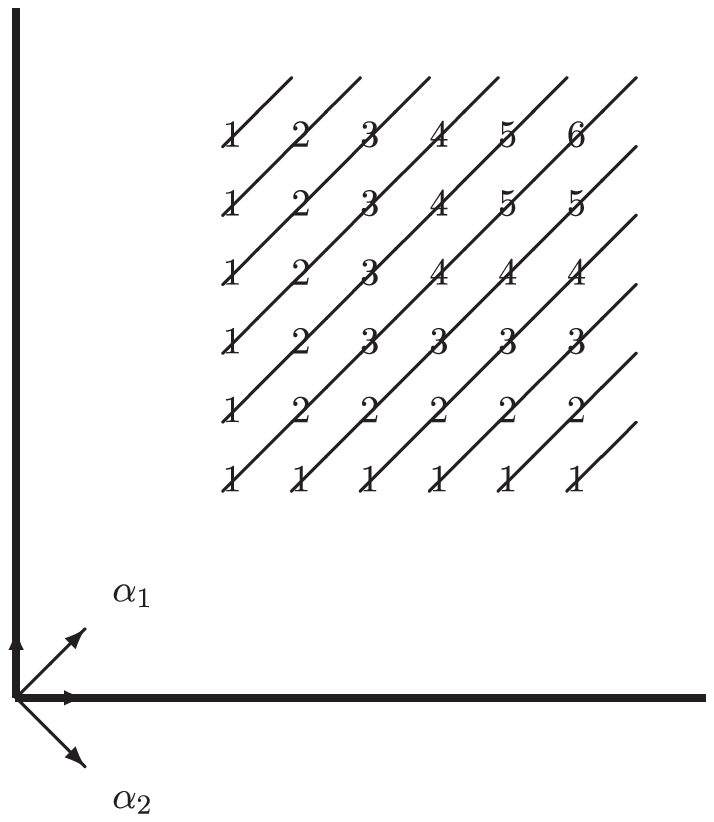}}
\caption{}\label{fig5}
\end{minipage}
\qquad
\begin{minipage}[b]{75mm}
\centerline{\includegraphics{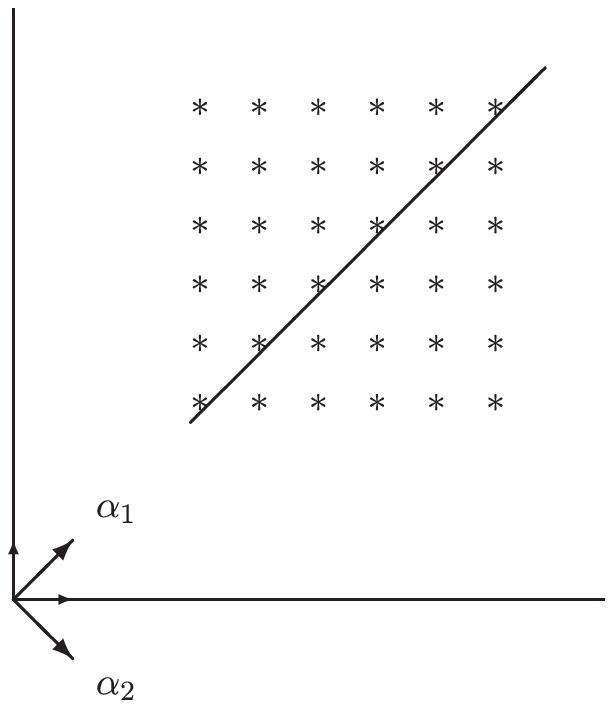}}
\caption{}\label{fig6}
\end{minipage}
\end{figure}

\noindent {\bf 6.2.}   These observations allow us to better
understand the restriction of $A_\bq$ to the pseudo dual pair
$H_1$=$Sp(2,\bR)$, $H_1'$= $GL(2,\bC)$ in $SL(4,\bR)$.
 The minimal $K$-type $V_K$ of $A_\bq$ has highest
weight~3$\alpha_1$ and is also irreducible under $K \cap H_1=
K\cap H_1'$. The $(\bh_1',K\cap H_1')$-submodule  generated by~$V_K$ is the minimal $H_1'$ type and is isomorphic to a
spherical principal series representation. We draw a diagram of
its $K \cap H_1'$-types using the same conventions as in the
previous sections. The $K \cap H_1'$-types are on  the black line
in Fig.~\ref{fig6}.

\begin{figure}[t]
\centering
\begin{minipage}[b]{75mm}
\centerline{\includegraphics{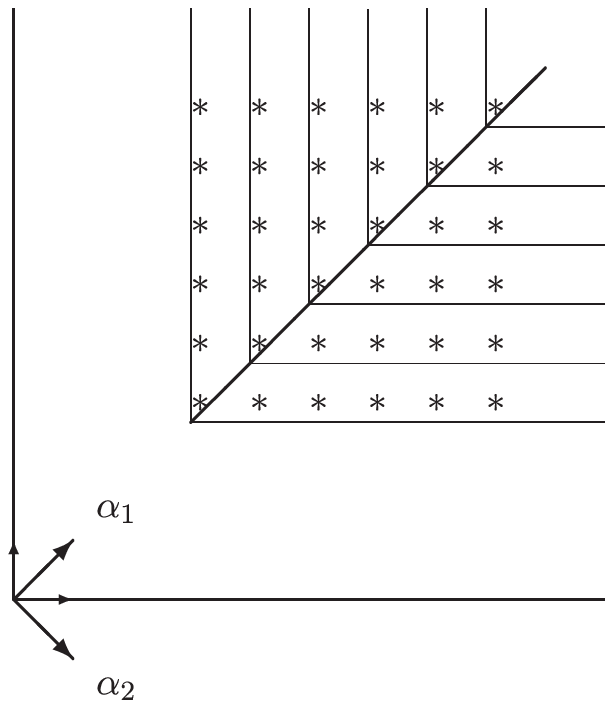}}
\caption{}\label{fig7}
\end{minipage}
\qquad
\begin{minipage}[b]{75mm}
\centerline{\includegraphics{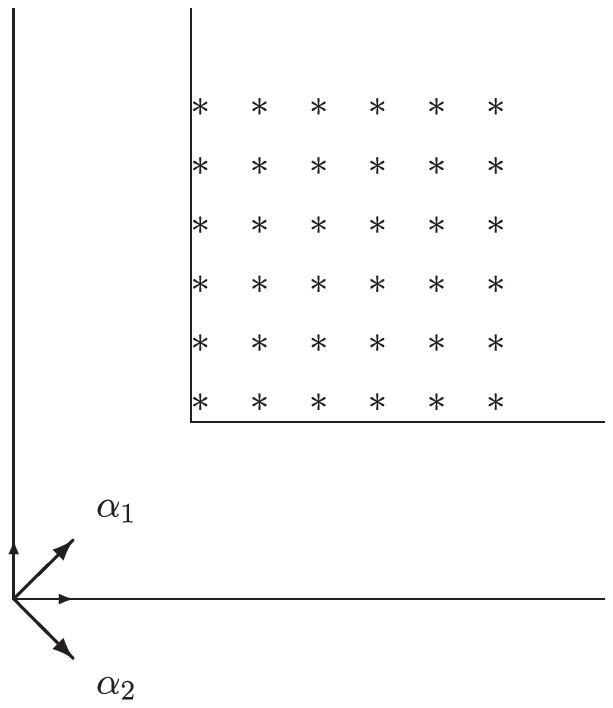}}
\caption{}\label{fig8}
\end{minipage}
\end{figure}

 Each $K\cap H_1=K\cap H_1'$-type of this representation is the minimal $K\cap
H_1$-type of an irreducible $(\bh _1, K \cap H_1)$-module as
indicated in Fig.~\ref{fig7}.

The $K\cap H_1$ types of the minimal $(\bh_1,K\cap H_1)$-module
generated by the minimal $K$-type have multiplicity one and are
indicated by the dots in Fig.~\ref{fig8}.

\begin{figure}[t]
\centerline{\includegraphics{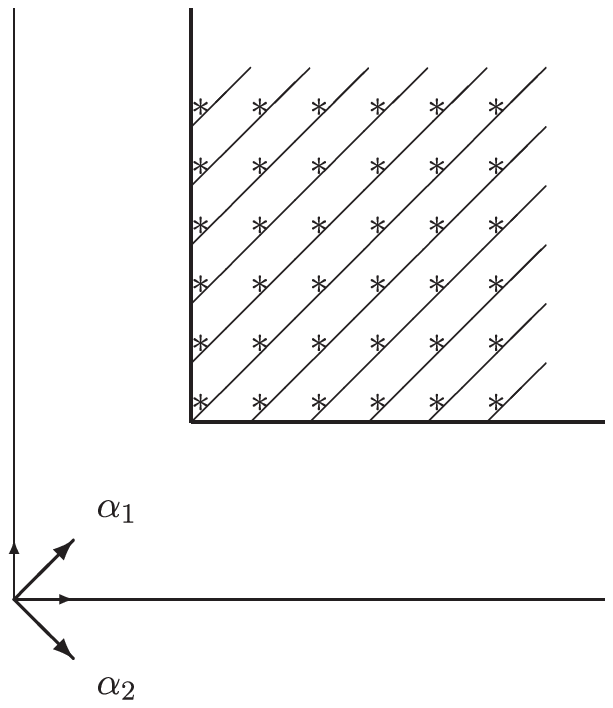}}
\caption{}\label{fig9}
\end{figure}

Each of the $K \cap H_1$-types of this $(\bh_1,K\cap
H_1)$-module is the minimal $K\cap H_1=K\cap H_1'$-type of a
$(\bh_1',K\cap H_1')$-module in the restriction of $A_\bq$ to
$H_1'$, as illustrated in Fig.~\ref{fig9}.

\section{More branching}\label{sec6}

In this section  we sketch the restriction to $H_1$ and $H_1'$ of
a representation $A_\bq(\lambda_1) $ with a parame\-ter~$\lambda $
which is no longer in the weakly fair range and hence the
representation is no longer irreducible and has a composition
series of length~2. We indicate a procedure to use our previous
techniques to compute the restriction of both composition factors
of $A_\bq(\lambda_1)$ to $H_1$ and $H_1'$. Using the local
isomorphisms
\begin{gather*}
SL(4,\bR) \sim  SO(3,3), \quad
Sp(2,\bR) \sim  SO(2,3), \quad
SL(4,\bR) \cap GL(2,\bC)  \sim  SO(2) \times SO(1,3)
\end{gather*}
we obtain a dif\/ferent proof of a result by T.~Kobayashi and B.~\O{}rsted of the branching of the minimal representation $\pi$ of
$SO(3,3)$ to $SO(3,1)\times SO(2)$ and to $SO(3,2)$ \cite{k-o}.

\medskip
\noindent {\bf 7.1.}   We use the notation introduced in~{\bf 3.1}.
Using the conventions on page~586 in~\cite{K-V} we denote the
character of $L$ by $\lambda = m(e_1+e_2+e_3+e_4)$. With this
parametrization the representa\-tion~$A_\bq(\lambda)$  is
irreducible and unitary for $m>-3$ (see page~588 in~\cite{K-V})
and $A_\bq$ corresponds to the parameter $m= 0$. The representation
$A_\bq(\lambda)$ for $m=3$ is nonzero.

To simplify the notation we denote the representation
$A_\bq(\lambda ) $ for $\lambda= m(e_1+e_2+e_3+e_4)$ by~$A_\bq(m)$.

We consider now the representation $A_\bq(-3)$. This
representation is nonzero and has a trivial $K$-type. It is not
irreducible, but has~2 composition factors. One composition factor
is $A_\bq(-2)$. The other composition factor is a unitarily
induced representation $\pi_0$ from one dimensional representation
of a maximal parabolic subgroup. It is a ladder representation and
the highest weights of its $K$-types are multiples of~$e_1$.
 See~\cite{k-o} for $SO(3,3)$ or~\cite{b-s-s} for~$GL(4,\bR)$.

The restriction of $A_\bq(-3)$ to $H_1 \cap K$ are $m_1(e_1 +e_2)+
2m_2   e_1$, $m_1, m_2 \in \bN.$

\medskip

\noindent {\bf 7.2.}  The same arguments as in Sections~\ref{sec3} and \ref{sec4}. allows us to compute the restriction of $A_\bq(-3)$ to
$H_1$ and $H_1'$. In this case the minimal $H_1$ type is a
spherical representation $\pi_0^{H_1}$ and comparing the
multiplicities of $K\cap H_1$-types of $\pi_0^{H_1}$ and of
$\pi_0$ we deduce that the restriction of $\pi_0$ to $H^1$ is
irreducible and equal to $\pi_0^{H_1}$. A similar argument show
that the restriction to $H_1'$ is direct sum of principal series
representations of $GL(2, \bC$).

Using the local isomorphism $SL(4,\bR )$ and $SO(3,3)$ we obtain a
new proof of the branching of $\pi_0$ determined by B.~\O{}rsted and
T.~Kobayashi in~\cite{k-o}.

\section{A conjecture}\label{sec7}

\noindent {\bf 8.1.} The examples in the previous section and
the calculations in~\cite{ko3} support the following conjecture:
Let $H$ be the connected f\/ixpoint set of an involution $\sigma$.
We write again $\bg = \bh \oplus \bs$.
 Let $A_\bq
(\lambda)$ be a representation, which is $K\cap H$-admissible and
thus decomposes discretely, when restricted to $H$. Suppose $\bq$,
$\bq^h $ are def\/ined by  $x_o \in T^H$. Since $A_{\bq} = A_{\bp}$
if $\bq$ and $\bp$ are conjugate under the compact Weyl group
$W_K$ we use the following

\medskip \noindent
{\bf Def\/inition.} Let $y_o \in T^H$ and let $\bp$, $\bp ^H $ be
well aligned parabolic subalgebras def\/ined by $y_o.$ We call the
well aligned parabolic subalgebras $\bp$, $\bp ^H $ {\it related}
to $\bq$, $\bq ^H$, if $x_o$ and $y_o$ are conjugate by an element
in the compact Weyl group $W_K$ of $K$ with respect to $T.$

\medskip

If $x_o$ and $y_o$ are not conjugate by an element in the Weyl
group $W_{K\cap H}$ of $(K^H, T^H)$ then the parabolic
subalgebras $\bq ^H$, $\bp ^H$ of H are not conjugate in $H$ and
thus we have up to conjugacy  at most $W_K/W_{H\cap K}$ dif\/ferent
{\it pairs of well aligned pairs of $\theta$-stable invariant
parabolic subalgebras which are related to $\bq$, $\bq ^H$.} If $G
= SL(4,\bR)$, $H= H_1$ and $(\bq, \bq \cap \bh_1 )$ is the pair of
well aligned parabolic subalgebras def\/ined by $x_0 = Q_2$, there
there are at most~2  related pairs of well aligned parabolic
subalgebras.

We expect the following Blattner-type formula to hold for the
restriction to $H$:

\medskip

\noindent
{\bf Conjecture.}
 {\it There exists a pair $\bp$, $\bp ^H $ of well aligned $\theta$-stable parabolic
 subgroups related to~$\bq$,~$\bq^H$ so that every $H$-type $V $ of~$A_\bq$
 is of the form $A_{\bp^H}(\mu )$ for a character
  $\mu $ of ${L^H} $
 and that }
\begin{gather*} \dim {\rm Hom}_{\bh,K^H} (V ,A_\bp )=
\sum_i  \sum_j (-1)^{s-j} \dim {\rm Hom}_{L \cap H} (H_j(\bu
\cap \bh, V), S^i(\bu \cap \bs )\otimes \bC_{\lambda_H}).
 \end{gather*}

\noindent
{\bf Remark.} Some of the characters $\mu $ in this formula may
be out of the fair range as def\/ined in~\cite{K-V} and hence
reducible.

If $H$ is the maximal compact subgroup $K$, then $|W_K/W_{H\cap K}| =
1$, all related pairs of well aligned parabolic subalgebra are
conjugate to $\bq$, $\bq ^H$ and hence we get the usual Blattner
formula (5.108b on page 376 in~\cite{K-V}).

In the example discussed in this paper, $G=SL(4,\bR)$, the
representation $A_\bq$ and $H$ the symplectic group considered in
Section \ref{sec3},  all related pairs of well aligned parabolic
subalgebras are conjugate to $\bq$, $\bq_{\bh_1}$ and thus we obtain
the Blattner type formula in Corollary~\ref{corollary4.4}.

\section{An application to  automorphic representations}\label{sec8}

 We use here our results  to give dif\/ferent
constructions of some known automorphic representations of $Sp(2,\bR)$
and $GL(2,\bC)$. We f\/irst explain the ideas in {\bf 7.1} %VI.1
in a more general setting. Again we may consider restrictions, this time in
the obvious way of restricting functions on locally symmetric
spaces to locally symmetric subspaces.

\medskip \noindent
{\bf 9.1.} Assume f\/irst that $G$ is a
semisimple matrix group and $\Gamma $ an arithmetic subgroup, $H$~a~semi\-simple subgroup of $G$.  Then $\Gamma_H = \Gamma \cap H$ is
an arithmetic subgroup of $H$.
 Let $ V_\pi \subset L^2(G/\Gamma)$ be an irreducible
$(\bg, K)$-submodule of $L^2(G/\Gamma)$. If $f \in V_\pi$ then $f$
is a $C^\infty$-function and so we def\/ine $f_{H}$ as the
restriction of f to $H/\Gamma_{H}$.

\begin{mylemma}
The map \begin{eqnarray*} \mbox{\rm RES}_{H}:\quad V_\pi &\rightarrow &
C^\infty(H/\Gamma_{H}) \\ f &\rightarrow &f_{H}
\end{eqnarray*}
is an $(\bh, K\cap H)$-map.
\end{mylemma}

\begin{proof} Let $h_t = \exp(t X_H)$, $h_o \in H$. Then
\begin{gather*}
 \rho(X_H)f(h_o) = \frac{d}{dt}f(h_f^{-1}h_0)_{t=0}= \frac{d}{dt}f_{H}(h_t^{-1}h_0)_{t=0} =
\rho(X_H)f_{H}(h_o) .\tag*{\qed}
\end{gather*} \renewcommand{\qed}{}
\end{proof}

 Suppose that the irreducible unitary $(\bg, K)$-module $\pi$ is a
submodule of $L^2(G/\Gamma)$ and that its restriction to $H$ is a
direct sum of unitary irreducible representations.

\begin{myprop}
Under the above assumptions $\mbox{\rm RES}_{H}(\pi)$ is nonzero and
its image is contained in the automorphic functions  on
$H/\Gamma_H$.
\end{myprop}

\begin{proof} Let $f_{H}$ be a function in  $\mbox{RES}_{H}(\pi)$. Then by Section~\ref{sec1} it is $K\cap H$-f\/inite and we
may assume that it is an  eigenfunction of the center of $U(\bh)$.

Let $||g||^2 = {\rm tr}(g^* g)$. Since $\sup\limits_{g\in G}|f(g)|\ ||g||^{-r }<
\infty $, the same is true for $f_{H}$ and so $f_{H}$ is an
automorphic function on $H/ \Gamma _{H}$.

The functions in the $(\bg, K)$-module $\pi \subset
L^2(G/\Gamma)$ are eigenfunctions of the center of the enveloping
algebra $U(\bg)$ and are $K$-f\/inite, hence analytic. Thus if $f $
is a $K$-f\/inite function in $\pi \subset L^2(G/\Gamma)$ then there
exists $W\in U(\bg)$ so that $Wf(e)\not = 0$. Hence $\mbox{RES}_{H}(Wf) \not = 0.$
\end{proof}

Instead of restricting the automorphic function f to the orbit of
$e/\Gamma $ under $H$ we may also consider the restriction to an
orbit of $\gamma /\Gamma, $ for rational~$\gamma $. Since the
rational elements are dense at least one of the restrictions is
not zero.
 So following Oda we consider the restriction correspondence for functions on
$ G/\Gamma $ to functions on $ \prod _g    H/H\cap
 g\Gamma g^{-1}$.
For rational $g$  the intersection $\Gamma \cap g \Gamma g^{-1}$
contains an arithmetic group $\Gamma '$ and $\Gamma ' \cap H $ is
an arithmetic subgroup. For more details see for example page~55
in~\cite{clozel-venky}.

\medskip

\noindent
{\bf 9.2.}  Now we assume that $G= GL(4, \bR)$ and that
$\Gamma \subset GL(4,\bZ)$ is a congruence subgroup.
 The groups $\Gamma_1 = \Gamma \cap H_1$ and $\Gamma_1'= \Gamma
\cap H_1'$ are arithmetic subgroups of $Sp(4,\bR)$, respectively
$GL(2,\bC)$. Recall the def\/inition of the $(\bg, K)$-module
${\mathcal A}_\bq$ from {\bf 3.1}. It is a submodule of $L^2(Z
\backslash G/\Gamma)$ for~$\Gamma$ small enough where $Z$ the
connected component of the center of $GL(4,\bR)$. We will for the
remainder of this sections consider it as an automorphic
representation in the residual spectrum~\cite{speh}. Then
 $\mbox{RES}_{H_1}({\mathcal A}_\bq)$ and $\mbox{RES}_{H_1'}({\mathcal A}_\bq)$
are nonzero. Its discrete summands are  contained in the space of
automorphic forms.

\begin{mytheorem}
The discrete summands  of the two representations $\mbox{\rm RES}_{H_1}({\mathcal A}_\bq)$ respectively \linebreak $\mbox{\rm RES}_{H_1'}({\mathcal A}_\bq)$ are subrepresentations  of the
discrete spectrum of $L^2(H_1/\Gamma _{H_1})$, respectively\linebreak
$L^2(H_1'/\Gamma _{H_1}')$.
\end{mytheorem}

\begin{proof}  All the functions in ${\mathcal
A}_\bq$ decay rapidly at the cusps. Since the cusps of
$H_1/\Gamma_{H_1}$ are contained in the cusps of $G/\Gamma $ this
is true for the functions in $\mbox{RES}_{H_1}({\mathcal
A}_\bq)$. Thus they are also contained in the discrete  spectrum.
\end{proof}

For $Sp(2,\bR)$ the representations constructed in the previous
theorem were f\/irst described by H.~Kim (see \cite{kim} and
\cite{schwermer}). For $GL(2,\bC)$ we obtain the stronger result

\begin{mytheorem}
The representations in the discrete spectrum of
 $\mbox{\rm RES}_{H_1'}({\mathcal A}_\bq)$ are in the cuspidal
 spectrum of $L^2(H_1'/\Gamma _{H_1'})$.
\end{mytheorem}

\begin{proof}
  By {\bf 5.2} the representations in the discrete spectrum of
 $\mbox{RES}_{H_1'}({\mathcal A}_\bq)$
 are unitarily induced
principal series representations and so by a result of Wallach
they are in fact cuspidal representations.
\end{proof}

The embedding of $H_1' = GL(2,\bC)$ into $SL(4,\bR)$ is def\/ined as
follows: Write $g = A+iB$ with real matrices $A$, $B$. Then
\[ g \rightarrow \left( \begin{array}{cc} A&B\\
                                           -B&A
 \end{array} \right ). \]
 Thus $\Gamma _{H_1'}$ is isomorphic to a congruence subgroup of $ GL(2,
 Z[i])$.

Since all the representations in the discrete spectrum of the
restriction of $A_\bq $ do have nontrivial $(\bh,K^H)$-cohomology
with respect to some irreducible f\/inite dimensional nontrivial
representation~$F$ we obtained the well known result
\cite{harder, rohlfs}.

\begin{mycor}
There exists a congruence subgroup $\Gamma \subset GL(2, Z[i])$
and a finite dimensional non-trivial representation F of $GL(2,
\bC)$ so that
\[ H^i(\Gamma, F) \not = 0 \qquad \mbox{for} \quad  i = 1,2. \]
\end{mycor}

\subsection*{Acknowledgements}

 We would like to thank T.~Kobayashi for helpful
discussions during a visit of the second author at RIMS and for
the suggestion to also include the restriction to $GL(2,\bC)$. The
second author would also like to thank the University of Southern
Denmark in Odense for its hospitality during which part of the
research was completed.
BS partially
supported by NSF grant DMS-0070561.

\pdfbookmark[1]{References}{ref}
\LastPageEnding

\end{document}